\documentclass[11pt,oneside]{amsart}
\usepackage[margin=1in]{geometry}
\usepackage{amsmath,amssymb,amsthm}
\usepackage[dvipsnames]{xcolor}
\usepackage[bookmarks,colorlinks,breaklinks]{hyperref}
\hypersetup{
  linkcolor=Blue,
  citecolor=BrickRed,
  filecolor=Purple,
  urlcolor=NavyBlue
}
\newtheorem{theorem}{Theorem}[section]
\newtheorem{lemma}[theorem]{Lemma}
\newtheorem{proposition}[theorem]{Proposition}
\newtheorem{corollary}[theorem]{Corollary}

\newcommand{\smallbreakheading}[1]{%
  \par\medskip\noindent\textit{#1}\par\smallskip
}
\makeatletter
\def\@setcopyright{}
\def\serieslogo@{}
\makeatother

\begin{document}

   \author{Amin  Bahmanian}
   \address{Department of Mathematics,
  Illinois State University, Normal, IL USA 61790-4520}\title{Sharp Embeddings of Layer-Rainbow Cubes and Orthogonal Pairs}


\begin{abstract}
We determine the exact extension threshold for one-factorizations of complete
tripartite $3$-uniform hypergraphs.  Equivalently, for $n>m$, every order-$m$
layer-rainbow cube embeds in one of order $n$ if and only if $n\ge2m$.  We
then study injectively orthogonal pairs, in which no ordered pair of symbols
is repeated.  Such pairs exist in every positive order except $2$, and every
order-$m$ pair with $m\ge3$ embeds in every order $n\ge3m$.  At the sharp
order $2m$, embedding occurs exactly when the pair has a disjoint mate.  The
symbol-pair graph carries three distinguished one-factorizations whose triple
intersections recover the cube, yielding several sufficient conditions for
sharp embedding, including a finite-group construction.  Finally,
layer-rainbow cubes are proper $n^2$-colorings of the square of the Hamming
graph, and the three one-factorizations of an orthogonal pair recover the
Hamming geometry of the cells.
\end{abstract}

\maketitle

\section{Introduction}
Several inequivalent notions of a ``Latin cube'' occur in the literature, so
we begin with the object considered here.  Let $L$ be an $n\times n\times n$
array.  A \emph{layer} of $L$ is obtained by fixing one coordinate.  If $L$
is filled with $n^2$ symbols so that every layer contains each symbol exactly
once, then $L$ is called a \emph{layer-rainbow cube}.  In particular, every
symbol occurs exactly $n$ times in the whole cube, once in each layer of any
fixed coordinate direction.

If Hamming distance counts the coordinates in which two cells differ, then
the $n$ cells carrying a fixed symbol have pairwise distance $3$, and no larger
set of cells can have this property because their first coordinates must all be
distinct.  Thus each symbol class is a maximum distance-$3$ code in $[n]^3$,
and a layer-rainbow cube partitions $[n]^3$ into $n^2$ such classes.

Layer-rainbow cubes were studied by Kishen \cite{MR34743}.  They are closely related to orthogonal arrays,
which have applications in statistics, coding theory, and cryptography
\cite{MR1693498,MR3495977}.  Symmetric layer-rainbow cubes have also been
considered recently in \cite{MR4665304}.

It is useful to distinguish this notion from a \emph{layer-Latin cube}, an
$n\times n\times n$ array on $n$ symbols in which every layer is a Latin
square.  A substantial part of the literature on layer-Latin completion
concerns non-completability and obstructions; see, for example,
\cite{MR677575,MR2838020,MR2902643}.  General finite embedding theorems are
known as well.  Cruse proved that every partial layer-Latin cube of order $n$
embeds in one of order $16n^4$ \cite{MR345854}, and Potapov later obtained an
order-$n^3$ embedding bound and related the problem to maximum distance
separable codes \cite{MR4453977}.

There is an equivalent hypergraph formulation that will be useful throughout.
The $n^3$ cells are the edges of the complete tripartite $3$-uniform
hypergraph $K^3_{n,n,n}$, and the cells carrying any one symbol form a perfect
matching.  Thus a layer-rainbow cube is the same object as a one-factorization
of $K^3_{n,n,n}$, and the embedding problem asks when a given
one-factorization of $K^3_{m,m,m}$ can be extended to one of
$K^3_{n,n,n}$.  Lemma~\ref{onefacrainbowlem} records this equivalence
formally.

Throughout the paper the three coordinate directions are distinguished.  If
$L$ has order $m$ and $U$ has order $n$, we say that $L$ \emph{embeds} in
$U$ if, after relabeling the three coordinate sets and the symbol set, $L$
occurs as an $m\times m\times m$ subcube of $U$.

The corresponding two-dimensional statement for complete Latin squares is
classical: for $n>m$, every Latin square of order $m$ embeds in one of order
$n$ if and only if $n\ge2m$ \cite{MR3495977}.
Our first theorem gives the exact three-dimensional analogue for layer-rainbow
cubes.

\begin{theorem}\label{layerrainbowembthm}
For $n>m$, a layer-rainbow cube of order $m$ can be embedded in a
layer-rainbow cube of order $n$ if and only if $n\ge2m$.
\end{theorem}

The same threshold is natural in every dimension.  A $d$-dimensional
layer-rainbow hypercube of order $n$ has $d$ coordinate sets of size $n$ and
$n^{d-1}$ symbols, each appearing once in every coordinate layer.  We
conjecture that, for $d\ge2$ and $n>m$, every order-$m$ layer-rainbow
hypercube embeds in one of order $n$ if and only if $n\ge2m$.  The lower
bound is necessary: the $m^{d-1}$ old symbols need $n-m$ further occurrences,
all in the $(n-m)^d$ cells whose coordinates are new, so
$m^{d-1}(n-m)\le(n-m)^d$, and hence $n\ge2m$.  Equivalently, we conjecture
the same threshold for extending one-factorizations of complete $d$-partite
$d$-uniform hypergraphs.

For related embedding results for graph and hypergraph factorizations, see,
for example, \cite{MR1249714,MR3910877,MR4728465}.

Our approach follows the same broad idea as Hilton's outline squares method
\cite{MR920647}.  The old and new coordinates divide the containing cube into
eight regions: one $m\times m\times m$ region, three regions with two old
coordinate directions, three with one old coordinate direction, and the
opposite $(n-m)\times(n-m)\times(n-m)$ region.  We first determine how many
copies of each symbol belong in these regions, without deciding their exact
positions.  We then split the coordinate vertices one at a time so that the
required numbers are distributed among individual layers.

Orthogonality is one of the classical themes in the theory of Latin squares,
going back to Euler's thirty-six officers problem \cite{MR3495977}.  For
$m>1$, a given pair of orthogonal Latin squares of order $m$ can be embedded
in a pair of every order $n\ge3m$ \cite{HiltonRodgerWojciechowski1992}.  The
factor $3$ is necessary: in the complementary
$(n-m)\times(n-m)$ subsquare, the remaining occurrences of the old symbols
from the two squares occupy two disjoint sets of $m(n-m)$ cells.  Hence
$2m(n-m)\le(n-m)^2$, so $n\ge3m$.  Related sharp results on orthogonal Latin
squares with aligned subsquares were obtained by Heinrich and Zhu
\cite{MR837956}; see also \cite{Lindner1976,DonovanYazici2014} for embedding
results concerning orthogonal partial Latin squares.

Orthogonality for higher-dimensional Latin-type arrays has also been studied;
see, for example, Ethier, Mullen, Panario, Stevens, and Thomson
\cite{EthierMullenPanarioStevensThomson2012}.  In our setting, for $n>1$,
ordinary orthogonality is impossible, since there are $n^4$ possible ordered
pairs of symbols but only $n^3$ cells.  We therefore impose the strongest possible
injectivity condition: no ordered pair is repeated.  We call such a pair
\emph{injectively orthogonal}, and usually simply \emph{orthogonal}.

We determine the possible orders of orthogonal pairs, prove that every
given order-$m$ pair embeds in every order $n\ge3m$, and characterize
embedding at the smallest possible order $2m$ by the existence of a disjoint
mate.  Thus the universal $3m$ range from the two-dimensional orthogonal
embedding theorem reappears in this setting, although some pairs attain the
minimum possible containing order $2m$.  The ordered symbol pairs form an
$m$-regular bipartite symbol-pair graph carrying three distinguished
one-factorizations, one for each coordinate direction; their triple
intersections recover the cells exactly.  This gives several sufficient
conditions for embedding at order $2m$, including a finite-group construction.
We conjecture that every orthogonal pair of order $m\ge3$ embeds in an
orthogonal pair of every order $n\ge2m$.  The results below settle every
$n\ge3m$ and characterize the endpoint $n=2m$ by the existence of a disjoint
mate; thus only the interval $2m<n<3m$ remains open in general.

The structures introduced here also admit further algebraic and geometric
interpretations.  In the companion paper we develop connections with the
Hamming association scheme and maximum distance separable codes, interpret
the three coordinate one-factorizations as mutually orthogonal simplex
systems, and study group-theoretic and affine-geometric constructions of
layer-rainbow cubes and orthogonal pairs.

\section{Preliminaries}\label{s:prereq}
Throughout, hypergraphs are allowed to have parallel edges; a hypergraph with
no parallel edges will be called \emph{simple}.  Thus a hypergraph
$G=(V,E)$ has a finite vertex set $V$ and an edge multiset $E$.  Parallel
copies having the same vertex set have the same \emph{edge type}.  For an
edge type $e$, write $(e)_G$ for the number of its copies in $G$; when $G$ is
fixed, we simply write $(e)$.  We write $\mathbf 1_e(v)=1$ when $v\in e$ and $0$ otherwise.

A tripartite $3$-uniform hypergraph has three pairwise disjoint vertex
parts, and every edge contains one vertex from each part. When the three
parts are $A,B,C$ and every possible triple occurs once, we write
$K^3_{A,B,C}$, or $K^3_{|A|,|B|,|C|}$. We call $\{x,y,z\}$ an $xyz$-edge.
The degree $\deg_G(v)$ counts incident edges with multiplicity.

For a color set $\kappa$, a $\kappa$-edge-coloring is a partition of the edge
multiset into color classes $G(c)$, one for each $c\in\kappa$.  For an edge
type $e$, write $(e)_c$ for the number of copies of $e$ in color $c$.  A \emph{matching}
is a set of pairwise vertex-disjoint edges, and it is \emph{perfect} if it
meets every vertex.  A perfect matching is also called a \emph{one-factor}.
The coloring is a \emph{one-factorization} if every color class is a
one-factor.  When there is no ambiguity we abbreviate $\deg_G(v)$ and
$\deg_{G(c)}(v)$ by $\deg(v)$ and $\deg_c(v)$.

\begin{lemma}\label{onefacrainbowlem}
A one-factorization of $K^3_{n,n,n}$ is equivalent to a layer-rainbow cube of
order $n$.
\end{lemma}
\begin{proof}
We let the three vertex parts be
$X=\{x_1,\dots,x_n\}$, $Y=\{y_1,\dots,y_n\}$, and
$Z=\{z_1,\dots,z_n\}$.  Given a one-factorization with $n^2$ colors, we place
color $c$ in cell $(i,j,k)$ exactly when the edge $x_i y_j z_k$ has color
$c$.  In the layer obtained by fixing $i$, the $n^2$ edges incident with
$x_i$ have distinct colors; since there are exactly $n^2$ colors, that layer
contains each color once.  The same argument applies after fixing $j$ or
$k$, so the resulting array is layer-rainbow.

Conversely, we color $x_i y_j z_k$ by the symbol in cell $(i,j,k)$ of a
layer-rainbow cube.  For any symbol $c$, the layer condition says that exactly
one $c$-colored edge meets each vertex in $X$, and likewise each vertex in
$Y$ and $Z$.  Hence every color class is a one-factor, and the colors form a
one-factorization.
\end{proof}

For an integer $a$ and a real number $b$, we write $a\approx b$ when
$\lfloor b\rfloor\le a\le\lceil b\rceil$.  A family of sets is
\emph{laminar} if any two of its members are disjoint or one contains the
other.  We use the following lemma of Nash--Williams \cite{MR0916377}: if
$d$ is a positive integer and $\mathcal A,\mathcal B$ are two laminar
families of subsets of a finite set $S$, then some $U\subseteq S$ satisfies
\[
 |U\cap W|\approx |W|/d
 \qquad \text{for } W\in\mathcal A\cup\mathcal B.
\]

\section{The one-cube embedding theorem}\label{proofmainthm}
Throughout this section, we let
\[
 \kappa_1=\{1,\dots,m^2\},\qquad
 \kappa_2=\{m^2+1,\dots,n^2\},\qquad
 \kappa=\kappa_1\cup\kappa_2.
\] 
For the construction, temporarily merge the $m$ old vertices in coordinate
part $j$ into one vertex $x_j$, and the $n-m$ new vertices into one vertex
$y_j$; call this \emph{amalgamation}.  We color the resulting six-vertex
hypergraph with parallel edges and then split each representative back into the
vertices it stands for.

\smallbreakheading{Coloring after amalgamation}

\begin{lemma}[Balancing three edge types]\label{basislem0}
Suppose that $n\ge2m\ge4$, and let $K=\{1,\dots,n^2-m^2\}$.  Let $G$ be
the tripartite $3$-uniform hypergraph having exactly three edge types
\[
 e_1=y_1x_2x_3,\qquad e_2=x_1y_2x_3,\qquad
 e_3=x_1x_2y_3,
\]
each with multiplicity $m^2(n-m)$.  Then the copies of these three edges can
be colored with the colors in $K$ so that, for every $c\in K$,
\begin{equation}\label{hardone}
 (e_j)_c+(e_k)_c\le m\quad \text{for }1\le j<k\le3,
 \qquad
 \sum_{j=1}^3(e_j)_c\ge3m-n.
\end{equation}
\end{lemma}
\begin{proof}
We write
\[
 N=|K|=n^2-m^2,\qquad T=m^2(n-m),\qquad
 \mu=\frac{T}{N}=\frac{m^2}{m+n}.
\]
We let $q=\lfloor\mu\rfloor$ and write $T=qN+R$, where $0\le R<N$.
Since $n\ge2m$,
\[
 \lceil\mu\rceil\le \left\lceil\frac m3\right\rceil
 \le \left\lfloor\frac m2\right\rfloor.
\]
Also
\[
 3\mu-(3m-n)=\frac{n(n-2m)}{m+n}\ge0.
\]
We let
\[
 L=\max\{0,3m-n\},\qquad t=\max\{0,L-3q\}.
\]
The preceding inequality gives $3R\ge tN$.  Because $R<N$, we have
$t\in\{0,1,2\}$.

We choose three subsets $S_1,S_2,S_3$ of $K$, each of size $R$, so that every
color belongs to at least $t$ of them.  This is immediate if $t=0$.  If
$t=1$, then $N\le3R$; we partition $K$ into three sets of size at most $R$
and enlarge each to size $R$.  If $t=2$, then $3(N-R)\le N$; we choose
pairwise disjoint sets $C_1,C_2,C_3$ of size $N-R$ and let
$S_j=K\setminus C_j$.

For $c\in K$, we define
\[
 (e_j)_c=
 \begin{cases}
 q+1,&c\in S_j,\\
 q,&c\notin S_j.
 \end{cases}
\]
For each $j$ these multiplicities sum to $qN+R=T$, so all copies of $e_j$
are colored.  Each multiplicity is at most $\lceil\mu\rceil$, and hence any
two sum to at most $m$.  Finally,
\[
 \sum_{j=1}^3(e_j)_c
 \ge 3q+t\ge L\ge3m-n.
\]
Thus \eqref{hardone} holds for every color.
\end{proof}

\begin{lemma}[Coloring the six-vertex hypergraph]\label{basislem}
Suppose that $n \ge 2m \ge 4$. Let $G=(V,E)$ be a tripartite $3$-uniform hypergraph with $n^3$ edges, whose parts are $\{x_1,y_1\}, \{x_2,y_2\}, \{x_3,y_3\}$, such that every edge type $e$ has multiplicity
\[
(e)=\prod_{j=1}^3 \Big( m^{\mathbf 1_e(x_j)} (n-m)^{\mathbf 1_e(y_j)} \Big).
\]
There exists a $\kappa$-edge-coloring of $G$ such that, for every color
$c\in\kappa$ and $j\in\{1,2,3\}$,
\begin{align}\label{maincolorcond}
\deg_c(x_j)=m, \qquad \deg_c(y_j)=n-m.
\end{align}
\end{lemma}
\begin{proof}
We let $\mathcal C=\{(1,2,3),(2,3,1),(3,1,2)\}$.  We call an edge a $t$-edge if it contains exactly $t$ of $x_1,x_2,x_3$.
We first color the $3$-edges so that $(x_1x_2x_3)_c=m$ for every
$c\in\kappa_1$.  After relabeling the color set in Lemma~\ref{basislem0}
as $\kappa_2$,
color the three types of $2$-edges so that, for $c\in\kappa_2$,
\[
 (e_j)_c+(e_k)_c\le m\quad \text{for }1\le j<k\le3,
 \qquad \sum_{j=1}^3(e_j)_c\ge3m-n.
\]
For every color $c\in\kappa_2$ and $(j,k,\ell)\in\mathcal C$, the
pairwise bound from Lemma~\ref{basislem0} makes the following quantity
nonnegative.  We color the $1$-edges by letting
\begin{align} \label{easyone}
(x_j y_k y_\ell)_c= m - (x_j x_k y_\ell)_c - (x_j x_\ell y_k)_c.
\end{align}
This is possible since for $(j,k,\ell) \in\mathcal C$,
\[
\sum_{c\in\kappa_2} \Big( m - (x_j x_k y_\ell)_c - (x_j x_\ell y_k)_c \Big) = m|\kappa_2| - 2 m^2 (n-m) = (x_j y_k y_\ell).
\]
For every color $c\in\kappa_2$, the lower bound in
Lemma~\ref{basislem0} also gives
$ n-3m+\sum_{j=1}^3(e_j)_c\ge0$.  We may therefore color the
$0$-edges so that
\begin{align} \label{easiest}
(y_1y_2y_3)_c= \begin{cases}
    n-m & \text{for } c\in\kappa_1,\\
    n-3m+\displaystyle\sum\nolimits_{(j,k,\ell)\in\mathcal C}(x_jx_k y_\ell)_c& \text{for } c\in\kappa_2,
    \end{cases}
\end{align}
This is possible because
\[
\begin{aligned}
\sum_{c\in\kappa_1} (n-m)
&+ \sum_{c\in\kappa_2}
\Bigl( n - 3m + \sum_{(j,k,\ell)\in\mathcal C}
(x_j x_k y_\ell)_c \Bigr) \\
&= |\kappa_1|(n-m)
   + |\kappa_2|(n-3m)
   + 3m^2(n-m)
   = (y_1y_2y_3).
\end{aligned}
\]

For every color $c\in\kappa_1$, the only $c$-colored edges are the
$3$-edge and the $0$-edge, with multiplicities $m$ and $n-m$, respectively.
Hence $\deg_c(x_j)=m$ and $\deg_c(y_j)=n-m$ for every $j$.

For every color $c\in\kappa_2$ and $(j,k,\ell)\in\mathcal C$, by
\eqref{easyone} and \eqref{easiest} we have
\begin{align*}
\deg_c(x_j) &= (x_j y_k y_\ell)_c + (x_j x_k y_\ell)_c + (x_j x_\ell y_k)_c = m,\\[1mm]
\deg_c(y_j) &= (y_j x_k x_\ell)_c + (y_j y_k x_\ell)_c + (y_j y_\ell x_k)_c + (y_1y_2y_3)_c\\
&= (x_k x_\ell y_j)_c + \big(m - (x_\ell x_j y_k)_c - (x_\ell x_k y_j)_c\big) \\
&\quad + \big(m - (x_k x_j y_\ell)_c - (x_k x_\ell y_j)_c\big) \\
&\quad + \big(n - 3m + \sum_{(p,q,h)\in\mathcal C} (x_p x_q y_h)_c\big)\\
&= n-m.
\end{align*}
Thus the edge-coloring of $G$ satisfies \eqref{maincolorcond}.
\end{proof}

\smallbreakheading{Splitting the amalgamated vertices}

A \emph{detachment} reverses amalgamation: split one amalgamated vertex and
distribute its incident edge copies among the new vertices.  The next
proposition does this one vertex at a time while ensuring that the given old
cube can be restored at the end.

\begin{proposition}[Controlled splitting]\label{inductiverthm}
Let $n\ge2m\ge4$ and $6\le\ell\le3n$.  There exists an $n^2$-edge-colored,
$\ell$-vertex, tripartite $3$-uniform hypergraph $G$ with parts
$P_j\cup Q_j$ for $j\in\{1,2,3\}$.  For each such $j$, we have
$P_j\cap Q_j=\varnothing$, $x_j\in P_j$, $y_j\in Q_j$,
$|P_j|\le m$, and $|Q_j|\le n-m$.  Moreover:
\begin{enumerate}
\item for every color $c\in\kappa$ and $j\in\{1,2,3\}$,
\[
 \deg_c(x_j)=m-|P_j|+1,\qquad
 \deg_c(y_j)=n-m-|Q_j|+1,
\]
and every other vertex has degree $1$ in color $c$;
\item every edge type $e$ has multiplicity
\[
(e)=\prod_{j=1}^3
 (m-|P_j|+1)^{\mathbf 1_e(x_j)}
 (n-m-|Q_j|+1)^{\mathbf 1_e(y_j)};
\]
\item if $c\in\kappa_1$, then every edge of color $c$ lies either entirely
in $P_1\times P_2\times P_3$ or entirely in
$Q_1\times Q_2\times Q_3$.
\end{enumerate}
\end{proposition}

\begin{proof}
For $\ell=6$, we let $P_j=\{x_j\}$ and $Q_j=\{y_j\}$.  Lemma~\ref{basislem}
gives (1) and (2).  Its construction also gives (3): colors in $\kappa_1$
are used only on the $3$-edge $x_1x_2x_3$ and the $0$-edge $y_1y_2y_3$.

Assume $6\le\ell<3n$ and that $G$ satisfies the three conditions.  At least one of
the six sets $P_j,Q_j$ has not yet reached its target size.  Relabeling the
three coordinate directions if necessary, suppose first that $|P_1|<m$.  We let
\[
 d=m-|P_1|+1.
\]
By an \emph{incidence at $x_1$} we mean a pair consisting of $x_1$ and a
particular copy of an edge containing it.  We let $\mathbb H(x_1)$ be the set of
all such incidences.  For each color $c\in\kappa$, we let $\mathbb H_c$
contain the incidences whose edge has color $c$.  For each edge type $f$
containing $x_1$, we let $\mathbb H_f$ contain the incidences belonging to
copies of $f$.  The two families
\[
 \{\mathbb H_c:c\in\kappa\},\qquad
 \{\mathbb H_f:f\ni x_1\}
\]
are partitions of $\mathbb H(x_1)$ and hence are laminar.  By
Nash--Williams' lemma there is $Z\subseteq\mathbb H(x_1)$ such that
\[
 |Z\cap A|\approx |A|/d
\]
for every set $A$ in either family.

We add a new vertex $\alpha$ to $P_1$ and, for every incidence
$(x_1,e)\in Z$, move that incidence from $x_1$ to $\alpha$.  We denote the
resulting hypergraph by $G'$.  Since $|\mathbb H_c|=\deg_c(x_1)=d$, we obtain
\[
 \deg'_c(\alpha)=1,\qquad
 \deg'_c(x_1)=d-1=m-|P'_1|+1.
\]
All other colored degrees are unchanged, proving (1).

We now let $f=x_1uv$ be an edge type and let $e=\alpha uv$ be the new type
obtained from it.  By (2), $(f)$ is divisible by $d$, and hence
\[
 (e)_{G'}=\frac{(f)}d,
 \qquad
 (f)_{G'}=(f)-\frac{(f)}d.
\]
These are precisely the products required in (2) after $|P_1|$ is increased
by one.  Finally, moving a vertex inside the same part $P_1$ cannot change
whether an edge lies in $P_1\times P_2\times P_3$ or in
$Q_1\times Q_2\times Q_3$, so (3) is preserved.

If instead $|Q_1|<n-m$, the same argument is applied to $y_1$ with
$d=n-m-|Q_1|+1$.  Repeating this splitting step reaches every
$\ell\le3n$ and completes the induction.
\end{proof}

\smallbreakheading{Completion of the proof}

\begin{proof}[Proof of Theorem~\ref{layerrainbowembthm}]
\emph{Necessity.} Each old color has $n-m$ remaining edges, all of which lie
in the opposite $K^3_{n-m,n-m,n-m}$.  Hence
$m^2(n-m)\le(n-m)^3$, so $n\ge2m$.

\emph{Sufficiency.} First consider $m=1$.  On coordinate sets $\mathbb Z_n$, we
color each edge $(x,y,z)$ by $(y-x,z-x)\in\mathbb Z_n^2$.  Each color class
is a one-factor of $K^3_{n,n,n}$, and $(0,0,0)$ has color $(0,0)$; after
relabeling, this extends the unique order-$1$ cube.

We now assume $m\ge2$ and apply Proposition~\ref{inductiverthm} with $\ell=3n$.
Then $|P_j|=m$ and $|Q_j|=n-m$ for every $j$.  Condition (2) says that every
possible tripartite edge has multiplicity one, so $G\cong K^3_{n,n,n}$.
Condition (1) says that every color has degree one at every vertex; hence the
color classes form a one-factorization.

It remains to realize the given old factorization.
By condition (3), every old color $c\in\kappa_1$ uses no mixed $P$--$Q$ edge.
Since it has degree one at each vertex of $P_1\cup P_2\cup P_3$, its edges in
$P_1\times P_2\times P_3$ form a one-factor.  The $m^2$ old colors therefore
contribute $m^3$ edges to this block, which is the total number of edges in
$K^3_{P_1,P_2,P_3}$.  Hence the $\kappa_1$-colored edges there form a
one-factorization of $K^3_{P_1,P_2,P_3}$.  We replace this factorization,
color by color, with the given one.  No edge outside
$P_1\times P_2\times P_3$ is
changed, and every colored degree remains one.  The resulting
one-factorization of $K^3_{n,n,n}$ is therefore the required extension.
\end{proof}

\section{Pairs of cubes}\label{s:orthogonalpairs}

We now consider simultaneous embeddings of two cubes.  As discussed in the
introduction, for $n>1$ ordinary orthogonality is impossible here: there are
$n^4$ possible ordered pairs of symbols but only $n^3$ cells.  We therefore use the
strongest possible substitute, requiring that no ordered pair occur twice.

Let $P=(L,M)$ be a pair of layer-rainbow cubes of order $n$, with symbol set
$S$ for $L$ and symbol set $T$ for $M$.  We say that $P$ is
\emph{injectively orthogonal} if no ordered pair $(s,t)\in S\times T$ occurs in
more than one cell.  We will usually say simply that $P$ is \emph{orthogonal}.

There is also a graph-coloring interpretation.  Let $G_n$ be the graph whose
vertices are the cells $[n]^3$, with two cells adjacent when they agree in at
least one coordinate.  A layer-rainbow cube is exactly a proper coloring of
$G_n$ with $n^2$ colors.  Under this identification, two cubes are orthogonal
exactly when the corresponding proper colorings are orthogonal in the sense
that no ordered pair of colors occurs at two vertices.

For an orthogonal pair $P=(L,M)$ with symbol sets $S$ and $T$, define
its \emph{symbol-pair graph} $\Gamma(P)$ to be the bipartite graph with parts
$S$ and $T$, where $s\in S$ is adjacent to $t\in T$ if the pair $(s,t)$
occurs in a cell.  Each $s\in S$ occurs in exactly $n$ cells, and
orthogonality implies that the corresponding $n$ symbols of $T$ are
distinct.  Hence $s$ has degree $n$.  The same argument applies to vertices of
$T$, so $\Gamma(P)$ is $n$-regular with $n^2$ vertices in each part.

Let $P=(L,M)$ be an orthogonal pair of order $m$, and let
$P'=(L',M')$ be one of order $n$.  We say that $P$ embeds in $P'$ if, after
relabeling the three coordinate sets and the two symbol sets independently,
$L$ and $M$ occur on the same $m\times m\times m$ subcube of $L'$ and $M'$.

\smallbreakheading{Existence}

\begin{theorem}\label{thm:orthpairsexist}
Orthogonal pairs exist for every positive order except $2$.
\end{theorem}

\begin{proof}
For order $1$ there is only one cell, so the assertion is immediate.  Suppose
$n\ge3$.  We work in $\mathbb Z_n$ and first construct permutations
$\sigma,\tau$ such that
\begin{equation}\label{eq:injectivepair}
 x\longmapsto \bigl(\sigma(x)-x,\tau(x)-x\bigr)
\end{equation}
is injective.  If $n$ is odd, let $\sigma(x)=2x$ and $\tau(x)=x$.  Since $2$
is invertible modulo $n$, both maps are permutations, and the ordered pair in
\eqref{eq:injectivepair} is $(x,0)$.  Now let $n=2k\ge4$.  Every element of $\mathbb Z_{2k}$ can be
written uniquely as $a$ or $a+k$, where $0\le a<k$, and we let
\[
 \sigma(x)=-x,\qquad \tau(a)=2a,\qquad \tau(a+k)=2a+1.
\]
The map $\sigma$ is a permutation.  The first $k$ values of $\tau$ are the
even residues modulo $2k$, and the
remaining $k$ values are the odd residues, so $\tau$ is a permutation.  Suppose
two elements $x,x'$ give the same first coordinate in \eqref{eq:injectivepair}.
Then $-2x\equiv-2x'\pmod{2k}$, so $x\equiv x'\pmod{k}$.  Thus either $x=x'$,
or the two elements are $a$ and $a+k$ for some $0\le a<k$.  In the latter
case their second coordinates in \eqref{eq:injectivepair} are $a$ and $a+1-k$, which are
distinct modulo $2k$ because $k\ge2$.  Hence
\eqref{eq:injectivepair} is injective for every $n\ge3$.

We define
\[
 A(x,y,z)=(y-x,z-x),\qquad
 B(x,y,z)=(y-\sigma(x),z-\tau(x)).
\]
The array $A$ is layer-rainbow.  If $x$ is fixed, its two components vary
independently with $y$ and $z$.  If $y$ is fixed, the first component determines
$x$, and then the second determines $z$; the case with $z$ fixed is analogous.
The same argument shows that $B$ is layer-rainbow: when $y$ is fixed, the first
component determines $\sigma(x)$ and hence $x$, and then the second determines
$z$; when $z$ is fixed, the second component determines $\tau(x)$ and hence
$x$, and then the first determines $y$.

To check orthogonality, suppose the cell $(x,y,z)$ contains $(a,b)$ in
$A$ and $(c,d)$ in $B$.  Then $\sigma(x)-x=a-c$ and
$\tau(x)-x=b-d$, so the injectivity of \eqref{eq:injectivepair}
determines $x$ uniquely.  The value $(a,b)$ in $A$ then determines $y=x+a$
and $z=x+b$.  Hence the ordered pair $((a,b),(c,d))$ determines the cell
$(x,y,z)$ uniquely, so no ordered pair of symbols is repeated.

It remains to exclude order $2$.  The two occurrences of any fixed symbol must
differ in all three coordinates.  Hence if one occurrence is in $(x,y,z)$,
the other is forced to the unique cell that differs from it in every
coordinate.  Thus any ordered pair of symbols appearing in one cell appears
again in that unique cell, so no orthogonal pair of order $2$ exists.
\end{proof}

\smallbreakheading{Embedding every given pair at orders at least $3m$}

For $n\ge3m$, we use incomplete orthogonal Latin squares.

Let $X$ be a set of size $n$ and let $H\subseteq X$ have size $m$.  An
\emph{incomplete Latin square of order $n$ with hole $H$} is an $n\times n$
array indexed by $X$, empty exactly on $H\times H$, such that every row or
column indexed by $X\setminus H$ contains every symbol of $X$ once, while
every row or column indexed by $H$ contains every symbol of $X\setminus H$
once and no symbol of $H$.  Filling $H\times H$ with any Latin square on $H$
therefore gives a Latin square of order $n$.

Two such squares $U^\circ,V^\circ$ are \emph{orthogonal} if, outside
$H\times H$, placing them on the same cells produces every ordered pair in
$(X\times X)\setminus(H\times H)$ exactly once.

The exact existence theorem is known: two incomplete orthogonal Latin squares
of order $n$ with a common hole of size $m$ exist exactly when $n\ge3m$,
apart from the exceptional case $(n,m)=(6,1)$; see Theorem~1.2 of Colbourn
and Zhu \cite{ColbournZhu1995}.  In particular, the needed pair exists whenever
$m\ge3$ and $n\ge3m$.

\begin{theorem}\label{thm:3mpluspair}
Let $m\ge3$ and let $P=(L,M)$ be an orthogonal pair of layer-rainbow
cubes of order $m$.  Then $P$ embeds in an orthogonal pair of order $n$
for every $n\ge3m$.
\end{theorem}

\begin{proof}
We fix $n\ge3m$.  We shall construct an order-$n$ orthogonal pair whose
prescribed $m\times m\times m$ subcube can be replaced by the pair $P$.

Let $X$ be a set of size $n$ and let $H\subseteq X$ have size $m$.  Choose two
incomplete orthogonal Latin squares $U^\circ,V^\circ$ of order $n$ with common
hole $H$.  Fill their holes independently with arbitrary Latin squares on $H$,
obtaining Latin squares $U,V$ on $X$.  These fillings are temporary and will
later be replaced by the prescribed pair.

We define
\begin{equation}\label{eq:imolscubes}
 A(x,y,z)=(U(x,y),U(x,z)),\qquad
 B(x,y,z)=(V(x,y),V(x,z)).
\end{equation}
We regard the symbol sets of $A$ and $B$ as two disjoint copies of
$X\times X$.

Both arrays are layer-rainbow.  For $A$, if $x$ is fixed, the two components
vary independently through a row of $U$.  If $y$ is fixed, the first component
determines $x$ from column $y$, and then the second determines $z$ from row
$x$.  The case with $z$ fixed is analogous.  The same argument, with $V$ in
place of $U$, shows that $B$ is layer-rainbow.

We next show that, in every coordinate layer indexed by an element of $H$, a
cell of $A$ has both coordinates of its entry in $H$ if and only if that cell
lies in the $H\times H\times H$ subcube.  If the cell lies in that subcube,
then both coordinates lie in $H$ because the filled $H\times H$ subsquare of
$U$ is a Latin square on $H$.  Conversely, suppose the cell lies in such a
layer but outside the subcube.  If $x\in H$ is fixed, then $y\notin H$ or
$z\notin H$, so $U(x,y)\notin H$ or $U(x,z)\notin H$.  If $y\in H$ is fixed,
then either $x\notin H$, in which case $U(x,y)\notin H$, or $x\in H$ and
$z\notin H$, in which case $U(x,z)\notin H$.  The case with $z\in H$ fixed is
analogous.  The same statement holds for $B$.

For $a,b\in X$, let
\[
 Q(a,b)=(U(a,b),V(a,b)).
\]
Since $U^\circ$ and $V^\circ$ are orthogonal outside the common hole, $Q$ is
injective on $(X\times X)\setminus(H\times H)$.  Moreover,
\begin{equation}\label{eq:holepaircriterion}
 Q(a,b)\in H\times H\quad\Longleftrightarrow\quad(a,b)\in H\times H.
\end{equation}
Indeed, outside the hole the two incomplete squares realize exactly the ordered
pairs in $(X\times X)\setminus(H\times H)$, while inside the filled hole both
entries lie in $H$.

We now show that two distinct cells outside $H\times H\times H$ cannot receive
the same ordered pair of entries from $A$ and $B$.  Suppose $(x,y,z)$ and
$(x',y',z')$ lie outside $H\times H\times H$ and receive the same ordered
pair.  Since $(x,y,z)$ is outside the subcube, at least one of $(x,y)$ and
$(x,z)$ lies outside $H\times H$.  Suppose first that
$(x,y)\notin H\times H$.  Equality of the first coordinates of the entries in
$A$ and $B$ gives $Q(x,y)=Q(x',y')$.  By \eqref{eq:holepaircriterion},
$(x',y')\notin H\times H$ as well.  Since $Q$ is injective outside
$H\times H$, we obtain $(x,y)=(x',y')$, and hence $x=x'$.  Equality of the
second coordinates of the entries in $A$ now gives $U(x,z)=U(x,z')$, so
$z=z'$.  Thus the two cells are equal.  If instead $(x,z)\notin H\times H$,
the same argument follows after interchanging $y$ and $z$.

We also show that a cell outside $H\times H\times H$ cannot have its entry in
$A$ lying in $H\times H$ and its entry in $B$ also lying in $H\times H$.  If
this happened at $(x,y,z)$, then both $Q(x,y)$ and $Q(x,z)$ would lie in
$H\times H$.  By \eqref{eq:holepaircriterion}, both $(x,y)$ and $(x,z)$ would
lie in $H\times H$, and hence $x,y,z\in H$, contradicting that the cell lies
outside $H\times H\times H$.

We now identify each of the three coordinate sets of $P=(L,M)$ with $H$, the
symbol set of $L$ with the copy of $H\times H$ used by $A$, and the symbol set
of $M$ with the copy of $H\times H$ used by $B$.  Replace the
$H\times H\times H$ subcube of $A$ by $L$, and the corresponding subcube of
$B$ by $M$.  By the property proved above, in every affected coordinate layer
the entries from $H\times H$ occur exactly in the portion being replaced.
Since $L$ and $M$ are themselves layer-rainbow cubes on those symbol sets, the
resulting arrays remain layer-rainbow.

It remains to check orthogonality.  Inside $H\times H\times H$, no
ordered pair is repeated because $P$ is orthogonal.  Outside
$H\times H\times H$, no ordered pair is repeated by the argument above.
Finally, an ordered pair cannot occur once inside and once outside: every
inside cell uses entries from the two copies of $H\times H$, while no outside
cell can do so simultaneously.  Hence the resulting order-$n$ pair is
orthogonal and contains $P$.
\end{proof}

The one-cube embedding theorem gives the necessary lower bound $n\ge2m$ for
any embedding into a larger order, while Theorem~\ref{thm:3mpluspair} settles
every $n\ge3m$.  Thus only $2m\le n<3m$ remains unresolved.

\smallbreakheading{The minimum order $2m$}

Let $P=(L,M)$ be an orthogonal pair of order $m$ with symbol sets $S$
and $T$.  A \emph{disjoint mate} of $P$ is another order-$m$
orthogonal pair $P^*$ on the same symbol sets that uses no ordered
symbol pair used by $P$; equivalently,
\[
 E(\Gamma(P))\cap E(\Gamma(P^*))=\varnothing.
\]

\begin{theorem}\label{thm:2mcriterion}
The pair $P$ embeds in order $2m$ if and only if it has a disjoint mate.
\end{theorem}

\begin{proof}
Suppose first that $P$ embeds in an order-$2m$ pair.  Consider the first
containing cube and fix a symbol of its prescribed order-$m$ subcube.  In every
coordinate layer indexed by one of the prescribed coordinates, that symbol
already occurs once inside the prescribed subcube, so it cannot occur elsewhere
in that layer.  Its remaining $m$ occurrences must therefore lie in the
$m\times m\times m$ subcube whose three coordinates all lie outside the
prescribed coordinate sets.  The $m^2$ symbols of the prescribed cube
contribute $m^2\cdot m=m^3$ occurrences there, exactly the number of cells in
that subcube.  Hence those symbols occupy every cell.  Each coordinate layer
indexed by a coordinate outside the prescribed sets contains each of these
symbols exactly once, so the restriction to this subcube is an order-$m$
layer-rainbow cube on the original symbol set.

Applying the same argument to the second containing cube gives another order-$m$
layer-rainbow cube on its original symbol set.  Together the two restrictions
form an orthogonal pair $P^*=(L^*,M^*)$.  No ordered symbol pair used by
$P$ can occur in $P^*$, since that would repeat an ordered pair in the containing
orthogonal pair.  Thus $P^*$ is a disjoint mate.

Conversely, suppose that a disjoint mate $P^*$ exists.  Take any order-$2$
layer-rainbow cube $C$ and replace each of its eight cells by an
$m\times m\times m$ block.  For each of the four symbols $q$ of $C$, choose a
set $S_q$ of $m^2$ symbols for the first large cube, with the four sets $S_q$
pairwise disjoint.  Independently choose pairwise disjoint sets $T_q$ of size
$m^2$ for the second large cube.

Each symbol $q$ occurs in exactly two cells of $C$.  For each $q$, place a
relabeled copy of $P$ in one of the corresponding two blocks and a relabeled
copy of $P^*$ in the other, using symbol sets $S_q$ and $T_q$.  Choose one of
the copies of $P$ to retain the prescribed labels.

The resulting arrays are layer-rainbow.  Fix a coordinate layer of either
large array.  The corresponding coordinate layer of $C$ contains each symbol
$q$ exactly once.  In the block belonging to that occurrence of $q$, the
corresponding coordinate layer of the inserted order-$m$ cube contains every
symbol of $S_q$, or every symbol of $T_q$, exactly once.  Since the sets for
different $q$ are disjoint, every symbol of the large array occurs exactly once
in the chosen layer.

Finally, if an ordered pair $(s,t)$ occurs, then $s\in S_q$ and $t\in T_q$
for a unique $q$: every block uses $S_q$ and $T_q$ with the same index, and the
families of sets are pairwise disjoint.  The pair can therefore occur only in
the two blocks corresponding to the two cells of $C$ containing $q$.  Those
blocks contain $P$ and $P^*$, which use disjoint ordered symbol pairs.  Hence no
ordered pair repeats.  The resulting order-$2m$ pair is orthogonal and
contains the prescribed copy of $P$.
\end{proof}

\smallbreakheading{The symbol-pair graph and its three one-factorizations}

The three coordinate directions give three distinguished one-factorizations of $\Gamma(P)$.  For
$h\in\{1,2,3\}$ and $a\in[m]$, let $F^h_a$ be the set of edges of $\Gamma(P)$
corresponding to cells in the layer whose $h$th coordinate is $a$.

\begin{theorem}\label{thm:threefactorizations}
For every orthogonal pair $P$ of order $m$, the families
\[
 \mathcal F_h=\{F^h_1,\ldots,F^h_m\},\qquad h=1,2,3,
\]
are one-factorizations of $\Gamma(P)$, and
\[
 |F^1_a\cap F^2_b\cap F^3_c|=1
 \qquad\text{for all }a,b,c\in[m].
\]
Conversely, let $\Gamma$ be a simple bipartite graph with two parts of size
$m^2$.  If $\Gamma$ has three one-factorizations whose perfect matchings can be
labeled $F^h_a$ so that this identity holds, then they determine an
orthogonal pair of order $m$, up to relabeling the coordinate layers and
the two symbol sets.
\end{theorem}

\begin{proof}
Each cell gives one edge of $\Gamma(P)$, and orthogonality ensures that
different cells give different edges.  In a fixed coordinate layer, every
symbol of each cube occurs exactly once, so the corresponding $m^2$ edges form
a perfect matching.  The $m$ layers in one coordinate direction partition the
cells and hence the edges.  Therefore each $\mathcal F_h$ is a
one-factorization.  Three specified layers, one in each coordinate direction,
meet in exactly one cell, which gives the displayed triple-intersection
identity.

Conversely, for each $(a,b,c)\in[m]^3$, let $s_{abc}t_{abc}$ be the unique edge
in $F^1_a\cap F^2_b\cap F^3_c$, and place $s_{abc}$ and $t_{abc}$ in cell
$(a,b,c)$ of two arrays.  Fixing one coordinate, say the first at $a$, uses
exactly the edges of the perfect matching $F^1_a$, so that layer contains every
symbol exactly once in each array.  The other two coordinate directions are
analogous.  If two cells received the same ordered pair of symbols, they would
correspond to the same edge.  But an edge belongs to a unique perfect matching
in each of the three one-factorizations, so it determines a unique triple
$(a,b,c)$.  Hence the two cells would be equal.  The arrays therefore form an
orthogonal pair.
\end{proof}

The theorem also recovers the coordinate geometry from the three
one-factorizations.  Each edge of $\Gamma(P)$ lies in a unique member of each
$\mathcal F_h$, and so receives a unique label $(a,b,c)\in[m]^3$.  The
triple-intersection condition makes this labeling a bijection from
$E(\Gamma(P))$ to $[m]^3$.  Under this bijection, two edges lie in different
perfect matchings in exactly $r$ of the three one-factorizations if and only if
their labels differ in exactly $r$ coordinates.  Thus the three
one-factorizations recover the Hamming scheme $H(3,m)$, that is, the relations
on $[m]^3$ determined by Hamming distance.

Fix distinct $h,h'\in\{1,2,3\}$ and perfect matchings $F^h_a$ and
$F^{h'}_b$.  As the remaining coordinate ranges over $[m]$, the corresponding
triple intersections are $m$ disjoint singletons whose union is
$F^h_a\cap F^{h'}_b$.  Hence
\[
 |F^h_a\cap F^{h'}_b|=m.
\]
Thus every perfect matching in one of the three one-factorizations meets every
perfect matching in another in exactly $m$ edges.

For a bipartite graph $\Gamma$ with parts $S$ and $T$, let $N_\Gamma(u)$ denote
the neighborhood of $u$.  For $U\in\{S,T\}$, define a bipartite graph
$A_U(\Gamma)$ with two copies $U^-$ and $U^+$ of $U$ by
\[
 u^-v^+\in E(A_U(\Gamma))
 \quad\Longleftrightarrow\quad
 N_\Gamma(u)\cap N_\Gamma(v)=\varnothing.
\]

\begin{lemma}\label{lem:auxmatchingmate}
Let $P$ be an orthogonal pair of order $m$ with symbol sets $S$ and $T$.
If either $A_S(\Gamma(P))$ or $A_T(\Gamma(P))$ has a perfect matching, then
$P$ has a disjoint mate and hence embeds in order $2m$.
\end{lemma}

\begin{proof}
Suppose $A_T(\Gamma(P))$ has a perfect matching.  It gives a permutation
$\beta$ of $T$ such that
\[
 N_\Gamma(t)\cap N_\Gamma(\beta(t))=\varnothing
 \qquad\text{for every }t\in T.
\]
Relabel the symbols of the second cube by $t\mapsto\beta(t)$.  Relabeling
preserves the layer-rainbow and orthogonal properties.  An edge $st$ of
the original symbol-pair graph becomes $s\beta(t)$.  If $s\beta(t)$ were also
an edge of the original graph, then $s$ would belong to both
$N_\Gamma(t)$ and $N_\Gamma(\beta(t))$, a contradiction.  Thus the relabeled
pair is a disjoint mate.  The proof for $A_S(\Gamma(P))$ is the same after
interchanging the two symbol sets.
\end{proof}

We shall use the following elementary bound.  If $\Gamma$ is an $m$-regular
bipartite graph with parts $S$ and $T$, each of size $m^2$, then every vertex
of $A_U(\Gamma)$ has degree at least $m-1$ for $U\in\{S,T\}$.
Indeed, fix $u\in U$.  The vertex $u$ itself shares a neighbor with $u$, and
for each of the $m$ neighbors of $u$, at most $m-1$ other vertices of $U$ share
that neighbor.  Thus at most $1+m(m-1)=m^2-m+1$ vertices $v\in U$ satisfy
$N_\Gamma(u)\cap N_\Gamma(v)\ne\varnothing$.  At least $m-1$ vertices remain,
which proves the bound.

\begin{proposition}\label{prop:structuralmate}
Let $P$ be an orthogonal pair of order $m\ge3$ with symbol sets $S$ and
$T$.  Then $P$ embeds in order $2m$ if either
\begin{enumerate}
\item[(i)] $\Gamma(P)$ contains no cycle of length four; or
\item[(ii)] for $U=S$ or $U=T$, and for any $u,v\in U$, there is an automorphism of
$\Gamma(P)$ (that is, a bijection of its vertices that preserves adjacency)
that maps $u$ to $v$ while mapping $S$ to $S$ and $T$ to $T$.
\end{enumerate}
\end{proposition}

\begin{proof}
Suppose first that $\Gamma(P)$ has no $4$-cycle.  Fix $t\in T$.  If
$s,s'\in N_\Gamma(t)$ are distinct, then
$N_\Gamma(s)\setminus\{t\}$ and $N_\Gamma(s')\setminus\{t\}$ are disjoint;
otherwise some $t'\ne t$ would give the $4$-cycle
$t,s,t',s',t$.  Each of the $m$ sets has size $m-1$, so exactly
$1+m(m-1)$ vertices of $T$ share a neighbor with $t$.  Hence every vertex of
$A_T(\Gamma(P))$ has degree $m-1$.  This graph is regular of positive degree,
so it has a perfect matching.  Lemma~\ref{lem:auxmatchingmate} applies.

Now suppose the second condition holds for $U\in\{S,T\}$.  The stated
automorphisms preserve whether two vertices of $U$ have disjoint neighborhoods,
so all vertices in the copy $U^-$ have the same degree in $A_U(\Gamma(P))$.
Because the defining relation is symmetric in its two vertices, all vertices in
$U^+$ have that same degree.  By the bound above, the degree is at least
$m-1>0$.  Thus $A_U(\Gamma(P))$ is a regular bipartite graph of positive degree
and has a perfect matching.  Lemma~\ref{lem:auxmatchingmate} again applies.
\end{proof}

A particularly simple sufficient condition comes from a finite group.

\begin{proposition}\label{prop:cayleymate}
Let $m\ge3$, and let $G$ be a finite group of order $m^2$.  Suppose, after identifying both
symbol sets with $G$, that
\[
 E(\Gamma(P))=\{(g,gd):g\in G,\ d\in D\}
\]
for some $D\subseteq G$ with $|D|=m$.  Then $P$ has a disjoint mate and hence
embeds in order $2m$.
\end{proposition}

\begin{proof}
The set
\[
 D^{-1}D=\{d_1^{-1}d_2:d_1,d_2\in D\}
\]
has size at most $1+m(m-1)=m^2-m+1<m^2$.  Choose
$h\in G\setminus D^{-1}D$.  Then $D\cap Dh=\varnothing$.  Indeed, if
$d=d'h$ for some $d,d'\in D$, then $h=d'^{-1}d\in D^{-1}D$, a contradiction.

Relabel the symbols of the second cube by $t\mapsto th$.  An edge $(g,gd)$ is
then replaced by $(g,gdh)$, so $D$ is replaced by $Dh$.  Since
$D\cap Dh=\varnothing$, the original and relabeled symbol-pair graphs are
edge-disjoint.  Thus the relabeled pair is a disjoint mate, and
Theorem~\ref{thm:2mcriterion} gives an embedding in order $2m$.
\end{proof}

\begin{corollary}\label{cor:onepairsharp}
For every $m\ge3$, there exists an orthogonal pair of order $m$ that embeds
in an orthogonal pair of order $2m$.
\end{corollary}

\begin{proof}
The construction in Theorem~\ref{thm:orthpairsexist} is of the form in
Proposition~\ref{prop:cayleymate}.  Identify both symbol sets with the additive
group $G=\mathbb Z_m^2$.  Its symbol-pair graph has
\[
 E(\Gamma(P))=\{(g,g+d):g\in G,\ d\in D\},
 \qquad
 D=\{(x-\sigma(x),x-\tau(x)):x\in\mathbb Z_m\}.
\]
The set $D$ has size $m$ by \eqref{eq:injectivepair}, so
Proposition~\ref{prop:cayleymate} applies.
\end{proof}

\bibliographystyle{amsplain}
\bibliography{references_partI}

\end{document}